\documentclass[11pt]{article}

\usepackage{amsmath,amssymb,mathtools}

\newtheorem{theorem}{Theorem}
\newtheorem{as}{Theorem}[section]
\newtheorem{proposition}[as]{Proposition}
\newtheorem{lemma}[as]{Lemma}

\newtheorem{corollary}[as]{Corollary}

\newcommand{\qed}{\hspace*{\fill} \rule{7pt}{7pt} \vspace*{5mm}}
\newcommand{\Proof}{\noindent{\bf Proof.}\ \ }

\begin{document}

\title{Perfect sequence covering arrays}

\author{
Raphael Yuster
\thanks{Department of Mathematics, University of Haifa, Haifa
31905, Israel. Email: raphy@math.haifa.ac.il}
}

\date{}

\maketitle

\setcounter{page}{1}

\begin{abstract}

An {\em $(n,k)$ sequence covering array} is a set of permutations of $[n]$ such that each sequence of $k$
distinct elements of $[n]$ is a subsequence of at least one of the permutations.
An $(n,k)$ sequence covering array is {\em perfect} if there is a positive integer $\lambda$ such that
each sequence of $k$ distinct elements of $[n]$ is a subsequence of precisely $\lambda$ of the permutations.

While relatively close upper and lower bounds for the minimum size of a sequence covering array are known, this is not the case for
perfect sequence covering arrays. Here we present new nontrivial bounds for the latter.
In particular, for $k=3$ we obtain a linear lower bound and an almost linear upper bound.

\vspace*{3mm}
\noindent
{\bf AMS subject classifications:} 05B40, 05B30, 05B15, 05A05\\
{\bf Keywords:} covering array; sequence covering array; completely scrambling set of permutations; directed t-design

\end{abstract}

\section{Introduction}

Let $2 \le k \le n$ be positive integers. Let $S_n$ denote the set of permutations of $[n]=\{1,\ldots,n\}$ and
let $S_{n,k}$ denote the set of all sequences of $k$ distinct elements of $[n]$.
An {\em $(n,k)$ sequence covering array} denoted by ${\rm SCA}(n,k)$, is a set $X \subseteq S_n$ such that each $\kappa \in S_{n,k}$ is a subsequence
of some element of $X$. Naturally, one is interested in constructing an ${\rm SCA}(n,k)$ which is as small as possible.
Thus, let $f(n,k)$ denote the minimum size of an ${\rm SCA}(n,k)$.

Sequence covering arrays have been extensively studied, see \cite{CCHZ-2013,KHL+-2012} and the references therein which also provide some
important applications of sequence covering arrays to the area of event sequence testing. Observe first that $f(n,2)=2$ as can be
seen by taking any permutation and its reverse.
A-priori, for a constant $k$, it is not entirely obvious that $f(n,k)$ grows with $n$, as each permutation covers $\binom{n}{k}$ sequences
while $|S_{n,k}|=k!\binom{n}{k}$.
However, more is known. The first to provide nontrivial bounds for $f(n,k)$ was Spencer \cite{spencer-1972} and various improvements on the upper and lower bounds
were sequentially obtained by Ishigami \cite{ishigami-1995,ishigami-1996},
F\"uredi \cite{furedi-1996}, Radhakrishnan \cite{radhakrishnan-2003}, and Tarui
\cite{tarui-2008}. The (asymptotic) state of the art regarding $f(n,3)$ is the upper bound by Tarui \cite{tarui-2008}
and the lower bound of F\"uredi \cite{furedi-1996}:
\begin{equation}\label{e:1}
\frac{2}{\log e} \log n \le f(n,3) \le (1+o_n(1))2 \log n\;. \footnote{Unless stated otherwise, all logarithms are in base $2$.}
\end{equation}
We note that the limit $f(n,3)/\log n$ exists \cite{furedi-1996,tarui-2008}, but apparently its value is not known.
For general fixed $k$, the best asymptotic upper and lower bounds are that of Spencer
\cite{spencer-1972} and Radhakrishnan \cite{radhakrishnan-2003}, respectively:
\begin{equation}\label{e:2}
(1-o_n(1))\frac{(k-1)!}{\log e} \log n \le f(n,k) \le \frac{k}{\log (\frac{k!}{k!-1})} \log n\;.
\end{equation}
We see that (\ref{e:2}) provides logarithmic upper and lower bounds for $f(n,k)$, so the order of magnitude of $f(n,k)$ for fixed $k$, is known.

A natural design-theoretic question that arises when studying sequence covering arrays is that of perfectness.
Let $X$ be an ${\rm SCA}(n,k)$. We call $X$ {\em perfect} if there exists an integer $\lambda$ such that each $\kappa \in S_{n,k}$ is a subsequence
of precisely $\lambda$ elements of $X$. We call $\lambda$ the {\em multiplicity} and denote $(n,k)$ perfect sequence covering arrays by ${\rm PSCA}(n,k)$ allowing
them to be multisets.
In design-theoretic terms, a ${\rm PSCA}(n,k)$ with multiplicity $\lambda$ is a $k-(n,n,\lambda)$ directed design, see \cite{CD-2006}
for the chapter on directed designs by Bennett and Mahmoodi.
Notice that a ${\rm PSCA}(n,k)$ exists for every $2 \le k \le n$ since $S_n$ is such. Let, therefore, $g^*(n,k)$ denote the minimum size of a ${\rm PSCA}(n,k)$
and observe the trivial bounds $k! \le f(n,k) \le g^*(n,k) \le n!$.

An easy observation is that $g^*(n,k)$ is a multiple of $k!$. Indeed, if each $k$-sequence is covered precisely $\lambda$ times, then the size of the
corresponding ${\rm PSCA}(n,k)$ is $\lambda k!$ since each permutation covers precisely $\binom{n}{k}$ sequences and there are $k!\binom{n}{k}$ sequences to cover.
So, we define the integer $g(n,k)=g^*(n,k)/k!$. Stated otherwise, $g(n,k)$ is the smallest $\lambda$ such that a $k-(n,n,\lambda)$ directed design exists.
Observe the trivial bounds $1 \le f(n,k)/k! \le g(n,k) \le n!/k!$.
We will also use the simple bounds $g(n,k) \ge g(n-1,k)$ and $g(n,k) \ge g(n,k-1)/k$.
Indeed, the former can be seen by taking any ${\rm PSCA}(n,k)$ and removing element $n$ from each permutation while the latter can be seen by taking the union of $k$ repeated copies of any ${\rm PSCA}(n,k)$.

Determining when $g(n,k)=1$ or, equivalently, when $f(n,k)=k!$, is an open problem. While clearly $g(k,k)=1$ and $g(n,2)=1$
it is a result of Levenshtein \cite{levenshtein-1991} that $g(k+1,k)=1$. 
It is also known that $g(6,4)=1$ \cite{MV-1999} and it is conjectured that $g(n,k) = 1$ only if $n \le k+1$ except for $k=2,4$ \cite{levenshtein-1991,MV-1999}.
The conjecture is known to hold for $k=3,5,6$ and it is also known that $g(7,4) > 1$ \cite{MV-1999}.
For general $k \ge 3$, a result of Chee et al. \cite{CCHZ-2013} shows that $g(2k,k) > 1$.

In Section 3 we determine the first (hence presently the only) exact bound of $g(n,k)$ which is not $1$ as we prove that $g(5,3)=2$.
However, our first main result is a lower bound for $g(n,k)$ which is much larger than the logarithmic lower bound for $f(n,k)$.
\begin{theorem}\label{t:1}
For all $n \gg k$, $g(n,k) > n^{k/2-o_k(1)}$. Furthermore, if $k/2$ is a prime then for all $n \ge k$ we have
$$
g(n,k) \ge \frac{\binom{n}{k/2}-\binom{n}{k/2-1}}{k!}\;.
$$
\end{theorem}
Notice that $g(n,4) \ge n(n-3)/48$ so coupled with the fact that $g(n,k) \ge g(n,k-1)/k$ we obtain, for every fixed $k \ge 4$, a polynomial in $n$ lower
bound for $g(n,k)$ while $f(n,k)$ is only logarithmic in $n$. Yet, Theorem \ref{t:1} does not give valuable input for the smallest nontrivial case $k=3$.
This is done in the next theorem, where we prove that $g(n,3)$ is at least linear in $n$ and at most quasi-linear in $n$.
\begin{theorem}\label{t:2}
For all $n \ge 3$, $n/6 \le g(n,3) \le Cn (\log n)^{\log 7}$ for some absolute constant $C$.
\end{theorem}
We note that the $\log 7 < 2.81$ can slightly be improved to any value strictly larger than $\log 6$ at the price of increasing $C$, but we
cannot eliminate it completely.

In the next section we prove our general lower bound, Theorem \ref{t:1}. The case $k=3$ and the proof of Theorem \ref{t:2} appear in Section 3.
The final section lists some open problems.
\section{Lower bounds}

Here we prove Theorem \ref{t:1}.
Let $X$ be a multiset of elements of $S_n$ and let $t$ be a positive integer.
We define the binary incidence matrix $A=A_{X,t}$ as follows.
The rows of $A$ are indexed by the elements of $S_{n,t}$, (all sequences of $t$ distinct elements of $[n]$)
and the columns of $A$ are indexed by $X$. For $\sigma \in X$ and $\kappa \in S_{n,t}$ we have $A[\kappa,\sigma]=1$ if
$\kappa$ is a subsequence of $\sigma$. Otherwise, $A[\kappa,\sigma]=0$.
We trivially have ${\rm rank}(A_{X,t}) \le |X|$.

We will prove Theorem \ref{t:1} for even values of $k$ such that $k/2$ is a prime.
We will then show that the result for other $k$ follows as a consequence.
Suppose now that $X$ is a ${\rm PSCA}(n,k)$ with multiplicity $\lambda$. Let $t=k/2$ and consider $A=A_{X,t}$.
Thus, we have ${\rm rank}(A_{X,t}) \le |X| = k!\lambda$.
We next consider the matrix $B=B_{X,k} = AA^T$. So clearly, ${\rm rank}(B) \le {\rm rank}(A) \le k!\lambda$.
Our goal is to obtain a lower bound for ${\rm rank}(B)$ which will imply a lower bound for $\lambda$.
Consider for example the case of $n=5$, $k=4$, $t=2$ and $\lambda=1$ where
a corresponding ${\rm PSCA}(5,4)$, which is also a $4-(5,5,1)$ directed design proving that $g(5,4)=1$, is given in Figure \ref{f:1}.
Figure \ref{f:2} shows the line of $B$ which corresponds to the sequence $12$.
\begin{figure}[h!]
$$
\begin{array}{cccc}
12345 & 12543 & 51423 & 41523\\
13524 & 15342 & 14325 & 54132\\
52134 & 21453 & 24135 & 42513\\
23514 & 25341 & 52431 & 42315\\
53124 & 31452 & 43512 & 34125\\
32154 & 45321 & 32451 & 35421
\end{array}
$$
\caption{A construction of a ${\rm PSCA}(5,4)$ with $\lambda=1$ showing that $g(5,4) = 1$.}
\label{f:1}
\end{figure}

\begin{figure}[h!]
$$
\begin{array}{c||c|c|c|c|c|c|c|c|c|c|c|c|c|c|c|c|c|c|c|c|}
& 12 & 13 & 14 & 15 & 23 & 24 & 25 & 34 & 35 & 45 & 21 & 31 & 41 & 51 & 32 & 42 & 52 & 43 & 53 & 54\\
\hline
12 & 12 & 8 & 8 & 8 & 4 & 4 & 4 & 6 & 6 & 6 & 0 & 4 & 4 & 4 & 8 & 8 & 8 & 6 & 6 & 6
\end{array}
$$
\caption{The row of $B=AA^T$ corresponding to the sequence $12$ for $A=A_{X,2}$ where $X$ is the ${\rm PSCA}$ from Figure \ref{f:1}.}
\label{f:2}
\end{figure}

To obtain a lower bound for ${\rm rank}(B)$, let us look first more carefully at the case $k=4$ (so $t=2$) but for general $n$.
There are only a few options for the entries of $B[ab,cd]$ where $ab,cd \in S_{n,2}$. Indeed, if $\{a,b\} \cap \{c,d\} = \emptyset$,
then $B[ab,cd]=6\lambda$ since the number of elements of $S_{n,4}$ in which $a$ precedes $b$ and $c$ precedes $d$ is precisely $6$.
If $a=c,b=d$, then $B[ab,ab]=12\lambda$ since in precisely half of the permutations of $X$, $a$ precedes $b$.
Similarly one can immediately check that if $a=c,b \neq d$ then $B[ab,ad] = 8\lambda$, and so on. We have listed all possible configurations and their respective
values in Figure \ref{f:3}.
\begin{figure}[t]
$$
\begin{array}{c||c|c|c|c|c|c|c|}
& cd & ab & ad & ba & ca & bd & cb\\
\hline
ab & 6\lambda & 12\lambda & 8\lambda & 0 & 4\lambda & 4\lambda & 8\lambda
\end{array}
$$
\caption{The values of $B[ab,\kappa]$ for the various types of $\kappa \in S_{n,2}$. Here $a,b,c,d$ are distinct.}
\label{f:3}
\end{figure}
We notice that each entry of $B$ is (obviously) a multiple of $\lambda$.
Let $C$ be obtained from $B$ by dividing each element by the gcd of all the entries of $B$.
In particular, this gcd is $2\lambda$ and notice that ${\rm rank}(C)={\rm rank}(B)$.
For a prime $p$, let ${\rm rank}_p(C)$ denote the rank of $C$ over the field ${\mathbb F}_p$.
For the case $p=2$, we see that over ${\mathbb F}_2$, $C$ is now the binary matrix with $C[ab,cd]=1$
if and only if $\{a,b\} \cap \{c,d\} = \emptyset$.

Consider the sub-matrix $C'$ of $C$ consisting of all rows $ab$ such that $a < b$ and all columns $cd$ such that $c < d$.
Then, we can view $C'$ as a matrix whose rows and columns are indexed by the unordered pairs of $[n]$, so over ${\mathbb F}_2$,
$C'$ is a binary matrix with $C'[X,Y]=1$ if and only if $\{a,b\} \cap \{c,d\} = \emptyset$.
Thus, ${\rm rank}_2(C') \le {\rm rank}_2(C) \le {\rm rank}(C) = {\rm rank}(B)$.
But one can now observe that $C'$ is precisely the set inclusion matrix of pairs versus subsets of order $n-2$.

Set inclusion matrices have been introduced by Gottlieb \cite{gottlieb-1966} and have been
extensively studied. Wilson \cite{wilson-1990} determined the rank of set inclusion matrices over finite fields - 
we next state his theorem. For integers $1 \le t \le \min\{r,n-r\}$ let $W_{t,r,n}$ denote the following matrix.
Its rows are indexed by all $t$-subsets of $[n]$ and its columns by all $r$-subsets of $[n]$ and we have $W[T,R]=1$ if $T \subseteq R$
and $W[T,R]=0$ otherwise.
\begin{lemma}\label{l:1}[Wilson \cite{wilson-1990}]
Let $p$ be a prime. Then ${\rm rank}_p (W_{t,r,n})$ is
$$
\sum_{i \in D(r,t)} \binom{n}{i}-\binom{n}{i-1}
$$
where $D(r,t)$ is the set of all integers $i$ such that $0 \le i \le t$ and $\binom{r-i}{t-i} \neq 0 \bmod p$.
\end{lemma}
\begin{corollary}\label{c:1}
Let $p$ be a prime. Then  ${\rm rank}_p (W_{t,r,n})$ is at least $\binom{n}{t}-\binom{n}{t-1}$.
\end{corollary}

So, in our case above, $C'$ equals $W_{2,n-2,n}$ (recall that the column indices are unordered pairs of $[n]$ but we can just rename them
by their complements, which are $(n-2)$-subsets of $[n]$). So by Corollary \ref{c:1} we obtain that $rank_2(C') \ge n(n-3)/2$.
It follows that ${\rm rank}(B) \ge n(n-3)/2$. Recalling also that ${\rm rank}(B) \le 24\lambda$ we have
$\lambda \ge n(n-3)/48$. Hence, $g(n,4) \ge n(n-3)/48$.

We now generalize the argument to all even $k \ge 6$ such that $t=k/2$ is a prime (thus an odd prime).
Consider $B[\kappa,\sigma]$ where $\kappa,\sigma \in S_{n,t}$.
Suppose first that $\kappa \cap \sigma = \emptyset$ (meaning that no element appears in both sequences).
Then the overall number of elements of $S_{n,k}$ that contain both $\kappa$ and $\sigma$ as subsequences is
precisely $\binom{k}{t}$, thus $B[\kappa,\sigma]=\lambda \binom{k}{t}$.
But notice that since $t=k/2$ is an odd prime, then $\binom{k}{t}$ is not divisible by $t$,
so $B[\kappa,\sigma]$ is not divisible by $\lambda t$.

Suppose next that $\kappa \cap \sigma \neq \emptyset$. Let $U = \kappa \cup \sigma$ be the set of symbols used in at least one of them
and notice that $t \le |U| \le k-1$.
Let $Q$ be the set of permutations of $U$ that is consistent with both $\kappa$ and $\sigma$, so $q \in Q$ if both $\kappa$ and $\sigma$
are subsequences of $q$. For example, suppose $k=6$, $\kappa = 123$ and $\sigma = 269$, then
$Q=\{12369,12639,12693\}$. Fix some $S \subset [n]$ with $U \cap S = \emptyset$ and $|U \cup S| = k$.
So, in the last example we can take, say, $S = \{4\}$.
Let $P$ be the set of permutations of $U \cup S$ that is consistent with both $\kappa$ and $\sigma$.
So, each element of $P$ is obtained by taking some $q \in Q$ and placing the elements of $S$ in some locations.
We therefore have that $|P| = |Q|\binom{k}{s}s!$ where $s=|S| \ge 1$ and that $B[\kappa,\sigma]=\lambda|P|$.
But now notice that $B[\kappa,\sigma]=\lambda|P|$ is divisible by $\lambda k$, hence by $\lambda t$.

We have shown that $B[\kappa,\sigma]$ is not divisible by $\lambda t$ if and only if $\kappa \cap \sigma = \emptyset$.
Let $C$ be obtained from $B$ by dividing each element by the gcd of all the entries of $B$
(and recall that this gcd is divisible by $\lambda$), so ${\rm rank}(C)={\rm rank}(B)$.
We see that over ${\mathbb F}_t$, $C$ is now a matrix with $C[\kappa,\sigma] \neq 0$ if and only if $\kappa \cap \sigma = \emptyset$
and furthermore, all nonzero entries of $C$ are equal to the same nonzero element of ${\mathbb F}_t$, call it $d$.
Consider the sub-matrix $C'$ of $C$ consisting of all rows $\kappa$ corresponding to increasing sequences and all columns $\sigma$ 
corresponding to increasing sequences.
Then, we can view $C'$ as a matrix whose rows and columns are indexed by the unordered $t$-subsets of $[n]$, so over ${\mathbb F}_t$,
$C'$ is a matrix with $C'[X,Y]=d$ if and only if $X \cap Y = \emptyset$.
Thus, ${\rm rank}_t(C') \le {\rm rank}_t(C) \le {\rm rank}(C) = {\rm rank}(B)$.
But now, $d^{-1}C'$ is the set inclusion matrix of $t$-subsets versus $n-t$ subsets, namely $W_{t,n-t,n}$.
So, by Corollary \ref{c:1}, ${\rm rank}_t(C') \ge \binom{n}{t}-\binom{n}{t-1}$.
It follows that ${\rm rank}(B) \ge \binom{n}{t}-\binom{n}{t-1}$.
Recalling also that ${\rm rank}(B) \le k!\lambda$ we have 
$\lambda \ge (\binom{n}{t}-\binom{n}{t-1})/k!$. Hence, $g(n,k) \ge (\binom{n}{k/2}-\binom{n}{k/2-1})/k!$.

We have thus proved that for all even $k$ such that $k/2$ is a prime and for all $n \ge k$, the statement in Theorem \ref{t:1} holds.
To end the theorem we just recall that $g(n,k) \ge g(n,k-1)/k$ since a ${\rm PSCA}(n,k)$ with multiplicity $\lambda$ is also a
${\rm PSCA}(n,k-1)$ with multiplicity $\lambda k$ and recall the fact that the primes are dense in the sense that for every integer $k \ge 2$
there is always a prime between $k$ and $k-O(k^{21/40})=k-o(k)$ \cite{BHP-2001}.
Hence we conclude that for all $n$ sufficiently large, $g(n,k) > n^{k/2-o_k(1)}$. \qed

\section{$g(n,3)$}

The following three lemmas prove Theorem \ref{t:2}.
\begin{lemma}\label{l:31}
$g(n,3) \ge n/6$.
\end{lemma}
\Proof
Suppose that $X$ is a ${\rm PSCA}(n,3)$, let $A=A_{X,2}$ be the incidence matrix of ordered pairs w.r.t. $X$ as defined in the previous section and let
$B=AA^T$. Since $A$ has $|X|=6\lambda$ columns, ${\rm rank}(B) \le 6\lambda$. As we cannot determine all elements of $B$, we will settle for a sub-matrix
of $B$ for which we can. Let $C$ be the sub-matrix of $B$ corresponding to the rows and columns indexed by the ordered pairs
$(i,n)$ for $i=1,\ldots,n-1$ and also by the ordered pair $(n,1)$, which will be the index of the last row and column.
(note: there are larger sub-matrices of $B$ with the property that all of their elements can be determined, but they do not yield larger rank).
So, $C$ is an $n \times n$ matrix. We will prove that $C$ is non-singular.

We observe that each diagonal entry of $C$ is $3\lambda$ since there are precisely $3\lambda$ elements of $S$ in which $i$ precedes $n$ for $i=1,\ldots,n-1$
and similarly there are $3\lambda$ elements of $S$ in which $n$ precedes $1$.
Similarly, $C[(i,n),(j,n)]=2\lambda$ for $i \neq j$ where $1 \le i,j \le n-1$,
$C[(1,n),(n,1)]=C[(n,1),(1,n)]=0$ and $C[(i,n),(n,1)]=C[(n,1),(i,n)]=\lambda$ for $i=2,\ldots,n-1$.
For simplicity, we divide all entries by $\lambda$ and set $C^* \coloneqq C/\lambda$. Figure \ref{f:4} is an example of $C^*$ in the case $n=6$.
It is not difficult to see by the matrix determinant lemma that ${\rm det}(C^*)=3(n+1)$ so ${\rm rank}(C)=n$, proving that ${\rm rank}(B) \ge n$ and that $\lambda \ge n/6$. \qed
\begin{figure}[t]
$$
\left(  
\begin{array}{rrrrrr} 
3 & 2 & 2 & 2 & 2 & 0\\
2 & 3 & 2 & 2 & 2 & 1\\
2 & 2 & 3 & 2 & 2 & 1\\
2 & 2 & 2 & 3 & 2 & 1\\
2 & 2 & 2 & 2 & 3 & 1\\
0 & 1 & 1 & 1 & 1 & 3
\end{array}
\right)
$$
\caption{The matrix $C^*$ for $n=6$.}
\label{f:4}
\end{figure}

\begin{lemma}\label{l:32}
Set $\lambda_1=1$ and $\lambda_r = 2(3^{\lceil r/2 \rceil}+1)\lambda_{\lceil r/2 \rceil}$ if $r \ge 2$.
Then, for $r \ge 1$ we have $g(3^r,3) \le \lambda_r$.
\end{lemma}
\Proof
We prove the lemma by induction on $r$ where the case $r=1$ holds since $g(3,3)=1$.
Notice that since $g(n,3)$ is monotone non-decreasing in $n$, we only need to prove $g(3^r,3) \le \lambda_r$ for even $r$,
since $\lambda_r=\lambda_{r+1}$ when $r$ is odd.
So, let $r$ be even and assume that for $n=3^{r/2}$ there is a ${\rm PSCA}(n,3)$ of multiplicity $\lambda_{r/2}$.
We will prove that there is a ${\rm PSCA}(n^2,3)={\rm PSCA}(3^r,3)$ of multiplicity $\lambda_r$.

Suppose $X$ is a ${\rm PSCA}(n,3)$ with multiplicity $\lambda=\lambda_{r/2}$ (hence $|X|=6\lambda$). We will construct 
a ${\rm PSCA}(n^2,3)$, denoted by $Y$, such that $|Y|=2(n+1)|X|$, and hence the lemma will follow by the definition of $\lambda_r$.

Our basic building block is a finite affine plane of order $n$, which exists since $n$ is a prime power.
This means, in particular, that there are $n+1$ partitions $P_1,\ldots,P_{n+1}$ of $[n^2]$, such that each $P_i$ consists of
$n$ parts of size $n$ each, denoted by $P_{i,j}$ for $j=1,\ldots,n$ and such that for any pair of distinct elements of $[n^2]$,
there is exactly one partition $P_i$ that contains both of them in the same part of $P_i$.

We construct $Y$ as a union of two sets $W,Z$ of $S_{n^2}$, where $|W|=|Z|=(n+1)|X|$.
We describe $W$ and then describe $Z$. $W$ will further be the union of $n+1$ sets $W_1,\ldots,W_{n+1}$
with $|W_i|=|X|$. We construct $W_i$ using $P_i$ and $X$.
Each element of $W_i$ will correspond to some $\sigma \in X$ as follows.
For each $P_{i,j}$, fix some total order of its $n$ elements (for example, the monotone increasing order).
For $\sigma \in S_n$, let $\sigma(P_{i,j})$ be the permutation of $P_{i,j}$ corresponding to $\sigma$.
Formally, if the total order of $P_{i,j}$ is $a_1,\ldots,a_n$ then
$\sigma(P_{i,j})$ is the permutation $a_{\sigma(1)},\ldots,a_{\sigma(n)}$.
For $\sigma \in X$ let $\sigma(P_i)$ be the concatenation of
$\sigma(P_{i,\sigma(1)}),\ldots,\sigma(P_{i,\sigma(n)})$. We call each part of this concatenation a {\em block}, so there are $n$ blocks of size $n$ each.
We observe that $\sigma(P_i) \in S_{n^2}$ and set $W_i=\{\sigma(P_i)\,:\, \sigma \in X\}$.
Thus, $W = \bigcup_{i=1}^{n+1}W_i$ is a well-defined subset of $S_{n^2}$.

Next, define $Z$ to be following ``reverse'' of $W$. For a totally ordered set $T$, its reverse, denoted $rev(T)$ is the the total order which places
the last element first, the second to last element second, and so on.
Now for $\sigma \in X$ let $q(\sigma(P_i))$ be the concatenation of $rev(\sigma(P_{i,\sigma(1)})),\ldots,rev(\sigma(P_{i,\sigma(n)}))$.
Set $Z_i=\{q(\sigma(P_i))\,:\, \sigma \in X\}$
and $Z = \bigcup_{i=1}^{n+1}Z_i$. Finally, let $Y = W \cup Z$ and observe that indeed $|Y|=2(n+1)|X|$ and $Y \subset S_{n^2}$.

To visualize our construction, consider for example the case $n=3$ with $X=S_3$ being the trivial ${\rm PSCA}(3,3)$ (with $\lambda=1$).
We will use the affine space of order $3$ formed of $P_1=\{123,456,789\}$, $P_2=\{147,258,369\}$, $P_3=\{159,267,348\}$, $P_4=\{168,249,357\}$.
Assume that in this listings, $P_{i,1}$ appears first, then $P_{i,2}$, then $P_{i,3}$ and that the listed order of each $P_{i,j}$ is the fixed total
order (we have used here the monotone increasing order). So, for example, for $\sigma=231 \in X$, we have, say $\sigma(P_{4,2})=492$
and $\sigma(P_4)$ is the concatenation of $\sigma(P_{4,2}),\sigma(P_{4,3}),\sigma(P_{4,1})$ so it is $492573681$.
Similarly, $q(\sigma(P_i))$ is $294 375 186$.

It remains to prove that each element of $S_{n^2,3}$ appears as a subsequence of precisely $2(n+1)\lambda$ elements of $Y$,
thereby proving that $Y$ is a ${\rm PSCA}(n^2,3)$ of multiplicity $\lambda_r$.
So, let $abc \in S_{n^2,3}$. We will distinguish between two cases. Assume first that $\{a,b,c\}$ is contained in some $P_{i,j}$
(in the case $n=3$ this means that $\{a,b,c\}$ is the whole $P_{i,j}$ but for larger $n$ this is strict containment).
Then, since $X$ is a ${\rm PSCA}(n,3)$ with multiplicity $\lambda$, we have that $abc$ appears precisely $\lambda$ times in $W_i$.
If $i' \neq i$, then $a,b,c$ appear in distinct blocks of each element of $W_{i'}$ (here we used the property of the affine plane).
So, again, since $X$ is a ${\rm PSCA}(n,3)$ with multiplicity $\lambda$, we have that $abc$ appears precisely $\lambda$ times in $W_{i'}$.
The exact same arguments apply for $Z_i$ and the $Z_{i'}$. Overall, $abc$ appears as a subsequence of precisely $2(n+1)\lambda$ elements of $Y$.

Assume next that $\{a,b,c\}$ is not a subset of any $P_{i,j}$.
Let $\gamma$ be the unique index such that $\{a,b\}$ is a subset of some part of $P_{\gamma}$,
let $\beta$ be the unique index such that $\{a,c\}$ is a subset of some part of $P_{\beta}$
and let $\alpha$ be the unique index such that $\{b,c\}$ is a subset of some part of $P_{\alpha}$.
Note that $\alpha,\beta,\gamma$ are indeed unique and distinct as follows from the properties of an affine plane.
As in the previous case, we have that if $i \notin \{\alpha,\beta,\gamma\}$ then $a,b,c$ appear in distinct blocks of each element of $W_{i}$
so we have that $abc$ appears precisely $\lambda$ times in $W_{i}$, and similarly for $Z_{i}$.
So $abc$ appears $2(n-2)\lambda$ times in $\bigcup_{i \in [n+1] \setminus \{\alpha,\beta,\gamma\}}(W_i \cup
Z_i)$. How many times does $abc$ appear as a subsequence in $W_\beta$? The answer is $0$, since in each element of $W_\beta$,
$a$ and $c$ appear in the same block while $b$ appears in another block. The same holds for $Z_\beta$.
How many times does $abc$ appear as a subsequence in $W_\alpha \cup Z_\alpha$? Since $bc$ are in the same block
of each element of $W_\alpha \cup Z_\alpha$ and since in precisely half of the elements of each of $W_\alpha$ and  $Z_\alpha$,
the block containing $a$ appears before the block containing both $b,c$ (we use here the fact that a PSCA of triples is trivially also
a PSCA of pairs), we have that precisely for half of the possible $\sigma$ precisely one of $\sigma(P_\alpha)$ or $q(\sigma(P_\alpha))$
contains $abc$ as a subsequence. So, overall, $abc$ appears as a subsequence in $W_\alpha \cup Z_\alpha$ precisely $|X|/2$ times.
The same argument holds for $W_\gamma \cup Z_\gamma$. In total, $abc$ is a subsequence of
$$
2(n-2)\lambda + 0 + |X| = 2(n+1)\lambda
$$
where the last equality follows from $|X|=6\lambda$.
We have thus proved that each $abc \in S_{n^2,3}$ is a subsequence of precisely $2(n+1)\lambda$ elements of $Y$, as required.
\qed

It is easy to prove by induction that for $r=2^t$ we have $\lambda_r = 2^{t-1}(3^r-1)$ hence for $n$ which is of the form
$3^{2^t}$ we obtain from Lemma \ref{l:32} that $g(n,3) \le \frac{1}{2}n\log_3 n$.
The next lemma provides an upper bound that applies to all values of $n$.
\begin{lemma}\label{l:33}
For all $n \ge 3$ we have $g(n,3) \le Cn (\log n)^{\log 7}$ for some absolute constant $C$.
\end{lemma}
\Proof
We first prove that the lemma holds for $n=3^r$ where $r$ is an integer.
Let $t=\lceil \log r \rceil$.
We will prove by induction that
$$
\lambda_r \le 7^t 3^r
$$
where $\lambda_r$ is as defined in Lemma \ref{l:32}.
Note that this holds for $\lambda_1=1$ and for $\lambda_2=8$.
Since $\lambda_r=\lambda_{r+1}$ when $r$ is odd, it suffices to prove $\lambda_r \le 7^t3^r$ when $r$ is odd.
Notice that if $r$ is odd, then $\lceil \log((r+1)/2) \rceil \le t-1$ so
by the definition of $\lambda_r$ we have for $r$ odd and the induction hypothesis that
\begin{eqnarray*}
\lambda_r & = & 2(3^{(r+1)/2}+1)\lambda_{(r+1)/2}\\
& \le & 2(3^{(r+1)/2}+1)7^{t-1} 3^{(r+1)/2}\\
& = & 2(3^{(r+1)/2}+1)(7/3)^{t-1}\cdot 3^{t-1} 3^{(r+1)/2}\\
& \le & (7/3) \cdot 3^{(r+1)/2} \cdot (7/3)^{t-1} \cdot 3^{t-1} \cdot 3^{(r+1)/2}\\
& = & 7^t 3^r\;.
\end{eqnarray*}
So, whenever $n$ is of the form $3^r$ we have that $g(n,3) \le 7^t 3^r$,
where $t=\lceil \log \log_3 n \rceil$.
If $n$ is not of this form,
let $n'$ be the unique power of $3$ such that $n \le n_0 < 3n$
and since
$g(n,3) \le g(n',3)$ we have $g(n,3) \le 7^t 3^r$ where $r=\log_3 n'$ and $t = \lceil \log \log_3 n' \rceil$.
Hence, $g(n,3) \le Cn (\log n)^{\log 7}$ for an absolute constant $C$.
\qed

We end this section with a proof that $g(5,3)=2$, which is currently the only explicitly determined value of $g(n,k)$ which is not one.
\begin{proposition}\label{p:1}
$g(5,3)=2$.
\end{proposition}
\Proof Recall from the introduction that $g(5,3) > 1$ \cite{MV-1999}, hence we only need to prove $g(5,3) \le 2$. We construct a ${\rm PSCA}(5,3)$ with $\lambda=2$.
It is not difficult to compute all sets of six permutations that cover a maximum number of sequences.
As it turns out, there are such sets that cover $56$ elements of $S_{5,3}$.
For example, the following is such:
$$
X=\{12345,43215,35214,14523,25413,53412\}
$$
The only sequences uncovered by $X$ are $132,231,154,451$. On the other hand, the sequences $123,321,145,541$ are each covered twice.
For $\sigma \in S_n$ and for $X \subseteq S_n$, let $X_{\sigma} = \{\pi\sigma\,:\, \pi \in X\}$.
Now, consider $\sigma=13254$. Then, for $X$ above we obtain that
$$
X_{\sigma} = \{13254,52314,24315,15432,34512,42513\}
$$
The only sequences uncovered by $X$ are $123,321,145,541$. On the other hand, the sequences $132,231,154,451$ are each covered twice.
Hence $X \cup X_{\sigma}$ is a ${\rm PSCA}(5,3)$.
\qed

\section{Open problems}

Theorem \ref{t:1} proves that for every fixed $k \ge 3$, $g(n,k)$ is lower bounded by a polynomial in $n$ whose exponent grows with $k$.
While it is not difficult to slightly improve upon the trivial upper bound $g(n,k) \le n!/k!$, it would be interesting to obtain
polynomial upper bounds for $g(n,k)$.

Theorem \ref{t:2} proves that $g(n,3)$ is at least linear and not more than quasi-linear in $n$. It would be interesting to determine the
right order of magnitude of $g(n,3)$.

Proving additional exact values of $g(n,k)$ which are not of unit multiplicity in addition to $g(5,3)$ also
seems challenging.

\section*{Acknowledgment}

The author thanks the referees for useful comments.

\bibliographystyle{plain}

\bibliography{references}

\begin{thebibliography}{10}

\bibitem{BHP-2001}
R.~Baker, G.~Harman, and J.~Pintz.
\newblock The difference between consecutive primes, {II}.
\newblock {\em Proceedings of the London Mathematical Society},
  83(03):532--562, 2001.

\bibitem{CCHZ-2013}
Y.~Chee, C.~Colbourn, D.~Horsley, and J.~Zhou.
\newblock Sequence covering arrays.
\newblock {\em SIAM Journal on Discrete Mathematics}, 27(4):1844--1861, 2013.

\bibitem{CD-2006}
C.~Colbourn and J.~Dinitz.
\newblock {\em Handbook of {C}ombinatorial {D}esigns}.
\newblock CRC press, second edition edition, 2006.

\bibitem{furedi-1996}
Z.~F{\"u}redi.
\newblock Scrambling permutations and entropy of hypergraphs.
\newblock {\em Random Structures {\&} Algorithms}, 8(2):97--104, 1996.

\bibitem{gottlieb-1966}
D.~Gottlieb.
\newblock A certain class of incidence matrices.
\newblock {\em Proceedings of the American Mathematical Society},
  17(6):1233--1237, 1966.

\bibitem{ishigami-1995}
Y.~Ishigami.
\newblock Containment problems in high-dimensional spaces.
\newblock {\em Graphs and Combinatorics}, 11(4):327--335, 1995.

\bibitem{ishigami-1996}
Y.~Ishigami.
\newblock An extremal problem of $d$ permutations containing every permutation
  of every $t$ elements.
\newblock {\em Discrete Mathematics}, 159(1-3):279--283, 1996.

\bibitem{KHL+-2012}
D.~Kuhn, J.~Higdon, J.~Lawrence, R.~Kacker, and Y.~Lei.
\newblock Combinatorial methods for event sequence testing.
\newblock In {\em Fifth International Conference on Software Testing,
  Verification and Validation (ICST)}, pages 601--609. IEEE, 2012.

\bibitem{levenshtein-1991}
V.~Levenshtein.
\newblock Perfect codes in the metric of deletions and insertions.
\newblock {\em Diskretnaya Matematika (English translation: Discrete
  Mathematics and Applications, 1992, 2:3, 241--258)}, 3(1):3--20, 1991.

\bibitem{MV-1999}
R.~Mathon and T.~Van~Trung.
\newblock Directed t-packings and directed t-{S}teiner systems.
\newblock {\em Designs, Codes and Cryptography}, 18(1-3):187--198, 1999.

\bibitem{radhakrishnan-2003}
J.~Radhakrishnan.
\newblock A note on scrambling permutations.
\newblock {\em Random Structures {\&} Algorithms}, 22(4):435--439, 2003.

\bibitem{spencer-1972}
J.~Spencer.
\newblock Minimal scrambling sets of simple orders.
\newblock {\em Acta Mathematica Hungarica}, 22(3-4):349--353, 1972.

\bibitem{tarui-2008}
J.~Tarui.
\newblock On the minimum number of completely 3-scrambling permutations.
\newblock {\em Discrete Mathematics}, 308(8):1350--1354, 2008.

\bibitem{wilson-1990}
R.~M. Wilson.
\newblock A diagonal form for the incidence matrices of $t$-subsets vs.
  $k$-subsets.
\newblock {\em European Journal of Combinatorics}, 11(6):609--615, 1990.

\end{thebibliography}

\end{document}